\documentclass[11pt,oneside]{amsart}
\usepackage{graphicx}
\usepackage{psfrag}
\usepackage{amscd}
\usepackage{epsfig}
\usepackage{amsmath,ifthen, amsfonts, amssymb, amsopn,amsthm}

\usepackage{color}

\theoremstyle{definition}

\begin{document}
\title{Non-Property (T) for SO(n,1)}
\author{Yongbin Zhou}

\begin{abstract}
We construct a concrete model for the measured wall of finite dinensional hyperbolic space, and construct the measure on it, so that the hyperbolic distance between two points equals to the measure of walls seperating them, up to a constant positive scalar. This offers a concrete proof for the fact that SO(n,1) is not T-group.
\end{abstract}

\maketitle

\section{Introduction}
For $n\ge2$, there is a $4n$-dimensional locally $CAT(-1)$ Riemannian manifold which is not homotopy equivalent to any $CAT(0)$ cube complex.~\cite{NR} This example is constructed by quaternionic hyperbolic space. One may ask if we can use hyperbolic space to construct examples of any dimension. His construction relies on the property (T) of (any lattice of) the isometry group of the quaternionic hyperbolic space. However, $SO(n,1)$, as isometry group of real hyperbolic space of dimension $n$, doesn't have property (T). So we can't just construct a example of any dimension in his way.\\
Previous constructions in some papers are abstract by some Haar measure on Lie group. Briefly speaking, they consider the isometry group of hyperbolic space, and get a quotient over stablizer of a codimension-1 totally geodesic space. We will embed the set of measured wall in euclidean spaces. We won't go into the median geometry or measured wall. One can find more about them in ~\cite{WALL}.

\section{Existence of the Measure}
We construct such a measure.
\subsection*{Hyperplanes and $dS^n$}
In this section, we introduce the $dS^n$ space to modelize the hyperplanes.
On \(\mathbb R^{n+1}\) let
\[
\langle x,y\rangle = -x_0y_0+x_1y_1+\cdots+x_ny_n .
\]
The hyperboloid model of hyperbolic \(n\)-space is
\[
I^n
=
\{x\in \mathbb R^{n+1}:\langle x,x\rangle=-1,\ x_0>0\}.
\]

Let
\[
\mathrm{dS}^n
=
\{u\in \mathbb R^{n+1}:\langle u,u\rangle=1\}
\]
be de Sitter space. For each \(u\in \mathrm{dS}^n\), define
\[
P_u
=
\{x\in I^n:\langle x,u\rangle=0\}.
\]
The orthogonal complement \(u^\perp\) has quadratic form signature \((n-1,1)\), so
\[
P_u=u^\perp\cap I^n
\]
is a totally geodesic subspace of \(I^n\) of dimension \(n-1\).

On the other hand, every totally geodesic hyperplane \(P\subset I^n\) is of the form
\[
P=I^n\cap V
\]
for a \(n\)-dimensional linear subspace \(V\subset \mathbb R^{n+1}\). Its quadratic form orthogonal complement \(V^\perp\) is one-dimensional, so there exists \(u\in\mathbb R^{n+1}\), unique up to a scalar, such that
\[
V^\perp=\mathbb R u, \qquad \langle u,u\rangle=1.
\]
Thus
\[
P=P_u.
\]

Therefore the space of totally geodesic (n-1)-dimensional hyperbolic hyperplanes in \(I^n\) is identified with \(\mathrm{dS}^n\), up to sign
\[
\Omega
\cong
\mathrm{dS}^n/\{\pm 1\}.
\]
The \(u\) and \(-u\) determine the same hyperplane.

We can parametrize de Sitter space by
\[
u=(\sinh r,(\cosh r)\,\omega),
\qquad
r\in \mathbb R,\quad \omega\in S^{n-1}.
\]
Easy verification:
\[
\langle u,u\rangle
=
-(\sinh r)^2+(\cosh r)^2|\omega|^2
=
-(\sinh r)^2+(\cosh r)^2
=
1.
\]
Where $|\omega|=1$, and $\sinh$ is a bijection from $\mathbb{R}$ to itself.
So
\[
\mathrm{dS}^n\cong \mathbb R\times S^{n-1}.
\]

\subsection*{G-Invariant Measure}
In this section, we show that the measure is G-invariant.
The parametrization
\[
u(r,\omega)
=
(\sinh r,\cosh r\,\omega),
\qquad
r\in\mathbb R,\quad \omega\in S^{n-1},
\]
results in
\[
\langle \partial_r u,\partial_r u\rangle=-1.
\]
For tangent vectors $v_1,v_2\in T_\omega S^{n-1}$,
\[
\langle d_\omega u(r,v_1),d_\omega u(r,v_2)\rangle
=
\cosh^2 r\, g_{S^{n-1}}(v_1,v_2).
\]
Where $g_{S^{n-1}}$ is the standard Riemannian metric tensor on $S^{n-1}$. We get the induced tensor
\[
g_{\mathrm{dS}}
=
-dr^2+\cosh^2 r\,g_{S^{n-1}}.
\]
This tensor is not a Riemannian metric tensor. However, what we concern is the measure. Take
\[
|d\operatorname{vol}_{\mathrm{dS}}|
=
\sqrt{|\det g_{\mathrm{dS}}|}\,dr\,d\omega.
\]
Because
\[
g_{\mathrm{dS}} =
\begin{pmatrix}
-1 & 0_{1\times(n-1)} \\
0_{(n-1)\times 1} & \cosh^2 r\, g_{S^{n-1}}
\end{pmatrix}.
\]
We get
\[
\det g_{\mathrm{dS}} = -\cosh^{2(n-1)}r\,\det g_{S^{n-1}},
\]
and as a result
\[
|d\operatorname{vol}_{\mathrm{dS}}| = \cosh^{n-1}r\,dr\,d\omega.
\]

Therefore
\[
\boxed{
d\mu(r,\omega)
=
\cosh^{n-1}r\,dr\,d\omega
}
\]

It is $SO(n,1)$-invariant on de Sitter space, equivalently on the space of totally geodesic hyperplanes in $I^n$.


\subsection*{Equal to Distance}
In this section, we show that the equality in the claim is true.
The measure is G-invariant, and $G$ act transitively, so the function
\[
F(x,y) = \mu(\{H\in \Omega | H \text{ separates } x \text{ and } y\})
\]
only relies on $d_{M}(x,y)$.
There exists a function
\[
f : [0,\infty) \to [0,\infty)
\]
such that
\[
F(x,y) = f(d_M(x,y)).
\]

Suppose that $z$ lies on the geodesic between $x$ and $y$. A hyperplane separates $x$ from $y$ if and only if it separates $x$ from $z$ or separates $z$ from $y$, and these two cases cannot hold at the same time, up to a set of measure zero. So,
\[
F(x,y) = F(x,z) + F(z,y).
\]

So
\[
f(a+b) = f(a) + f(b).
\]

It is obvious that $f$ is measurable and finite on any bounded interval, so it must be linear. So
\[
f(t) = ct
\]
for some constant $c=c(n)\ge0$. And $c>0$ is obvious: by taking $|r|$ small and $\omega$ ranging in a small area on $S^{n-1}$, we can separate $o=(1,\dots,0)$ and another point very far away.

\section{$\mathrm{SO}(n,1)$ Is Not Property (T) by Conditionally Negative Kernels}

Let
\[
G=\mathrm{SO}(n,1),
\qquad
M=I^n.
\]
Then $G$ is the isometry group of $M$, where $M$ is the hyperboloid model of hyperbolic space. Fix a basepoint $o\in M$. Define
\[
\psi(g):=d_{M}(o,go),
\]
and
\[
K(g,h)=\psi(g^{-1}h)=d_{M}(go,ho).
\]

We will show that $K$ is CNK(conditionally negative kernel) but unbounded, if we just put the same method on hyperbolic space and SO(n,1). Then $G$ does not have property $(T)$. ~\cite{NR} relied on the property $(T)$, so we can't construct example by hyperbolic space in the same way. The condition we will use from NR's paper requires that G is finitely generated, but I believe there is a way to fix it even if $\mathrm{SO}(n,1)$ is not finitely generated.

Recall the following things copied from ~\cite{NR}:
\[
\textbf{Definition}\quad
\text{A conditionally negative kernel on a set } V \text{ is a function }
f : V \times V \longrightarrow \mathbb{R}
\]
\[
\text{such that for any finite subset } \{v_1,\ldots,v_n\} \subset V
\text{ and any real numbers } \{\lambda_1,\ldots,\lambda_n\}\] \[\sum_i \lambda_i = 0
\]
\[
\text{the following inequality holds:}
\]
\[
\sum_{i,j} \lambda_i \lambda_j f(v_i,v_j) \leq 0.
\]
A conditionally negative kernel on a group $G$ is a conditionally negative kernel on the set of elements of $G$ such that for any $g,h,k$ in $G$,
$f(gh,gk)=f(h,k)$.

A finitely generated group $G$ has Kazhdan's property (T) (or is a $T$--\textbf{group}) if and only if every conditionally negative kernel on $G$ is bounded.~\cite{KA}

\subsection*{Step 1: Set CNK: $d_{M}$}
In this section we prove $d_{M}$ is set CNK.
Let $\Omega$ be the space of totally geodesic (n-1)-dimensional hyperbolic hyperplane in $M$.
We have constructed a $G$-invariant measure $\mu$ on $\Omega$ such that
\[
\mu\{H\in \Omega | H \text{ separates } x \text{ and } y\}
=
c\, d_{M}(x,y)
\]
for some constant $c>0$. Let's set $c=1$.

Define
\[
\Phi(x)(H)=1,\,iff\,H\,separates\,x\,and\,o,\,otherwise\,0
\]
so that $\Phi(x)\in L^2(\Omega,\mu)$, because the hyperbolic distance between $x$ and $o$ in $M$ is finite, and $\Phi(x)(\_)$ is defined almost everywhere on $\Omega$.  

Then almost everywhere on $\Omega$
\[
|\Phi(x)-\Phi(y)|=1,\,iff\,H\,separates\,x\,and\,y,\,otherwise\,\,0
\]
So $|\Phi(x)-\Phi(y)|\in L^2(\Omega,\mu)$, because the hyperbolic distance between $x$ and $y$ in $M$ is finite, and $|\Phi(x)-\Phi(y)|$ is defined almost everywhere on $\Omega$.
\[
\|\Phi(x)-\Phi(y)\|^2
=
\mu\{H |  H \text{ separates } x \text{ and } y\}
=
 d_{M}(x,y).
\]

On M
\[
d_{M}(x,y) = \|\Phi(x)-\Phi(y)\|^2.
\]

Now let $\lambda_1,\dots,\lambda_m\in\mathbb R$ with $\sum_i \lambda_i=0$.
Then
\[
\sum_{i,j}\lambda_i\lambda_j d(x_i,x_j)
=
\sum_{i,j}\lambda_i\lambda_j \|\Phi(x_i)-\Phi(x_j)\|^2.
\]

For Hilbert space
\[
\sum_{i,j}\lambda_i\lambda_j \|f_i-f_j\|^2
=
-2\left\|\sum_i \lambda_i f_i\right\|^2,
\]
we get
\[
\sum_{i,j}\lambda_i\lambda_j d(x_i,x_j)
=
-{2}\left\|\sum_i \lambda_i \Phi(x_i)\right\|^2
\le 0.
\]

Hence $d_{M}$ is a conditionally negative kernel.

\subsection*{Step 2: Group CNK: K}

In this section we prove $K$ is group CNK.
For $g_1,\dots,g_m\in G$ and $\sum_i \lambda_i=0$, by CNK on the set,
\[
\sum_{i,j}\lambda_i\lambda_j K(g_i,g_j)
=
\sum_{i,j}\lambda_i\lambda_j d_{M}(g_i o, g_j o)
\le 0.
\]
And for $g,h,k\in G$, G is the isometry group of M, so
\[
K(gh,gk)
=
d_{M}(gho,gko)
=
d_{M}(ho,ko)
=
K(h,k)
\]

Thus $K$ is conditionally negative kernel on group $G$.

\subsection*{Step 3: Unboundedness of K}

Walking alone a geodesic line through $o$ in M, we get a one-parameter set of points ${p_t}$ with $(a_t)_{t\in\mathbb R}\subset G$ be a one-parameter subgroup acting $o$ into $p_t$ with
\[
p_0=o
\]
Then
\[
d_{M}(o,a_t o)=|t|.
\]

Hence
\[
K(a_{t},e)=d_M(a_to,o)=|t|\to \infty \quad \text{as } t\to\infty,
\]
so $K$ is unbounded. Therefore SO(n,1) is not T-group.

\end{document}